\documentclass[a4paper,11pt]{article}

\usepackage[latin1]{inputenc}
\usepackage{amsfonts,amsmath,amssymb,amsthm,graphicx,epsfig,float}
\usepackage[T1]{fontenc}
\usepackage[english,francais]{babel}
\usepackage[all]{xy}

\usepackage{amsthm}
\newtheorem{Thm}{Théorème}
\newtheorem{Pro}{Proposition}[section]
\newtheorem{Lem}{Lemme}[section]

\theoremstyle{definition}

\theoremstyle{remark}

\newcommand{\Cc}{\mathbb{C}}
\newcommand{\Nn}{\mathbb{N}}

\newcommand{\Zz}{\mathbb{Z}}
\newcommand{\Pp}{\mathbb{P}}

\newcommand{\Hcal}{\mathcal{H}}

\title{{\bf Minoration de la dimension de Hausdorff du courant de Green}}
\author{Henry De Thélin}
\date{}

\begin{document}
\maketitle

\def\figurename{{Fig.}}%
\def\proofname{Preuve}% for AMS-\LaTeX
\def\contentsname{Sommaire}%

\begin{abstract}

Nous donnons une minoration de la dimension de Hausdorff du support du courant de Green associé à certaines applications méromorphes.

\end{abstract}

\selectlanguage{english}
\begin{center}
{\bf{ }}
\end{center}

\begin{abstract}

We give a lower bound for the Hausdorff dimension of the support of the Green current associated to some meromorphic maps.

\end{abstract}

\selectlanguage{francais}

Mots-clefs: dynamique complexe, courants, dimension de Hausdorff.

Classification: 32H50, 32U40, 37C45, 37D25.

\section*{{\bf Introduction}}
\par

A partir d'un endomorphisme holomorphe $f$ de $\Pp^2(\Cc)$ de degré $d \geq 2$, J.-E. Forn{\ae}ss et N. Sibony ont construit le courant de Green $T$ qui est limite au sens des courants de $\frac{(f^n)^* \omega}{d^n}$ où $\omega$ est la forme de Fubini-Study de $\Pp^2(\Cc)$ (voir \cite{FS1} et \cite{FS2}). J.-Y. Briend a montré dans \cite{Br} que ce courant a un potentiel höldérien. En particulier la dimension de Hausdorff de son support est strictement plus grande que $2$ (voir le corollaire 1.7.4 de l'article de N. Sibony \cite{S}).

Lorsque $f$ est méromorphe dominante de $\Pp^2(\Cc)$ de degré $d \geq 2$, le courant de Green $T$ existe toujours par les travaux de J.-E. Forn{\ae}ss et N. Sibony mais à priori son potentiel n'est pas continu. Un des objectifs de cet article va être de présenter des situations où l'on peut malgré cela minorer la dimension de Hausdorff de son support par une quantité strictement plus grande que $2$. 

Considérons tout d'abord $(X, \omega)$ une surface kählérienne compacte et $f: X \longrightarrow  X$ une application birationnelle. On note $d_1$ son premier degré dynamique et on suppose que $d_1 > 1$. Dans ce contexte les courants de Green $T^{\pm}$ de $f$ et $f^{-1}$ existent par J. Diller et C. Favre (voir \cite{DF}). Dans \cite{BD}, E. Bedford et J. Diller ont introduit une condition sur l'ensemble d'indétermination $I$ de $f$ qui permet de définir une mesure $\mu=T^+ \wedge T^-$ et ils ont montré qu'elle est hyperbolique. R. Dujardin a prouvé que son entropie est égale à $\log d_1$ (voir \cite{Du}). Elle est donc d'entropie maximale par T.-C. Dinh et N. Sibony (voir \cite{DS1} et \cite{DS2}). 

Nous allons montrer que sous cette condition de Bedford-Diller, le courant de Green a son support de dimension de Hausdorff strictement plus grande que $2$. Plus précisément, on a 

\begin{Thm}{\label{bedford-diller}}

Soit $(X, \omega)$ une surface kählérienne compacte et $f: X \longrightarrow  X$ une application birationnelle qui vérifie la condition de Bedford-Diller. On note $d_1$ son premier degré dynamique et on suppose que $d_1 >1$. Soit $T^{\pm}$ ses courants de Green et $\mu=T^+ \wedge T^-$ sa mesure d'entropie maximale. Par \cite{BD} cette mesure est hyperbolique et si $\chi_1 > 0 > \chi_2$ désignent ses exposants de Lyapounov, on a

$$\mbox{dim}_{\Hcal} ( \mbox{Supp} T^+) \geq 2 + \frac{\log d_1}{\chi_1} > 2.$$

\end{Thm}

Ce théorème provient d'un résultat plus général qui est valable en dimension quelconque. En effet, on va montrer le

\begin{Thm}{\label{general}}

Soit $(X, \omega)$ une variété kählérienne compacte de dimension $k$ et $f: X \rightarrow X$ une application méromorphe dominante. On suppose que les degrés dynamiques de $f$ vérifient $d_0 \leq \cdots < d_s > \cdots \geq d_k$ avec $1 \leq s \leq k-1$ et que l'on a une mesure $\mu$ d'entropie maximale $\log d_s$ avec $\int \log d(x,I) d \mu(x) > - \infty$ (en particulier $\mu$ est hyperbolique avec $s$ exposants strictement positifs et $k-s$ strictement négatifs par \cite{Det1} et \cite{Dup}).

Alors, si on a convergence exponentielle  de $\frac{(f^n)^* \omega^s}{d_s^n}$ vers un courant de Green $T^+$ (dans le sens où il existe $\alpha_2 > 0$ tel que pour toute $(k-s,k-s)$ forme lisse $\Phi$ on ait

$$\exists n_0 \mbox{  ,  } \forall n \geq n_0  \mbox{  ,  } \left| \langle \frac{(f^n)^* \omega^s}{d_s^n} - T^+, \Phi \rangle \right| \leq e^{- \alpha_2 n}),$$

on a 

$$\mbox{dim}_{\Hcal} ( \mbox{Supp} T^+) \geq 2(k-s) + \frac{\log d_s}{\chi_1}.$$

Ici $\chi_1$ désigne le plus grand exposant de Lyapounov de $\mu$.

\end{Thm}

Pour montrer le premier théorème à partir de celui-là, il s'agira de prouver que sous la condition de Bedford-Diller, la convergence de  $\frac{(f^n)^* \omega}{d^n}$ vers le courant de Green est exponentielle. En effet, si $f$ vérifie la condition de Bedford-Diller, on a bien $\int \log d(x,I) d \mu(x) > - \infty$ (voir \cite{BD}).

Le point crucial pour démontrer le théorème ci-dessus est la présence de la mesure hyperbolique $\mu$. En effet, nous avons vu dans \cite{DetNgu} que l'on peut approcher sa dynamique par des horseshoes. Grâce à cette propriété, on peut construire une lamination par des variétés stables de dimension $k-s$. Nous verrons que la convergence  exponentielle de $\frac{(f^n)^* \omega^s}{d_s^n}$ vers le courant de Green $T^+$ implique que le support de $T^+$ contient cette lamination. Enfin, on minorera la dimension de Hausdorff de cette lamination en utilisant la dynamique et en particulier les exposants de Lyapounov de $\mu$.

Précisons maintenant ces idées. Dans tout ce qui suit nous reprendrons les résultats et les notations de \cite{DetNgu} et en particulier la thérie de Pesin développée dans le paragraphe 1.1 de ce papier.

On fixe $X$ une variété complexe compacte de dimension $k$, $f: X \longrightarrow X$ une application méromorphe dominante et nous noterons $I$ son ensemble d'indétermination. 

Comme dans \cite{DetNgu}, on munit $X$ d'une famille de cartes $(\tau_x)_{x \in X}$ qui vérifient $\tau_x(0)=x$ et on considère $\mu$ une mesure de probabilité avec $\int \log d(x,I) d \mu (x)>-\infty$. On suppose $\mu$ invariante par $f$, ergodique et hyperbolique, c'est-à-dire que ses exposants de Lyapounov vérifient : 

$$\chi_1\geq  \cdots \geq \chi_{s}>0>\chi_{s+1}\geq \cdots \geq \chi_{k} \mbox{  pour un  } 1\leq s \leq k-1.$$ 

On note $\Omega=X \setminus \cup_{i \geq 0} f^{-i}(I)$. Comme $\mu$ ne charge pas $I$ et qu'elle est invariante par $f$, $\mu$ est une probabilité de $\Omega$.

On définit l'extension naturelle:

$$\widehat{\Omega}:= \{ \widehat{x}=( \cdots, x_0, \cdots , x_n , \cdots) \in \Omega^{\Zz} \mbox{ , } f(x_{n})=x_{n+1} \}.$$

Dans cet espace, $f$ induit le décalage à gauche $\widehat{f}$ et si on note $\pi$ la projection canonique $\pi(\widehat{x})=x_0$ on notera $\widehat{\mu}$ l'unique probabilité invariante par $\widehat{f}$ qui vérifie $\pi_{*} \widehat{\mu}=\mu$.

Dans toute la suite, on pose $f_x= \tau_{f(x)}^{-1} \circ f \circ \tau_x$ qui est définie au voisinage de $0$ quand $x$ n'est pas dans $I$. Comme dans \cite{DetNgu}, le cocycle auquel nous allons appliquer la théorie de Pesin est:
\begin{equation*}
\begin{split}
A :\ & \widehat{\Omega} \longrightarrow M_k(\Cc)\\
& \widehat{x} \longrightarrow Df_x(0)\\
\end{split}
\end{equation*}

où $M_k(\Cc)$ est l'ensemble des matrices carrées $k \times k$ à  coefficients dans $\Cc$ et $\pi(\widehat{x})=x$. Grâce à l'hypothèse d'intégrabilité de la fonction $\log d(x,I)$, on a un théorème du type Oseledets (voir le théorème $2$ et la proposition $1.1$ dans \cite{DetNgu}). Nous noterons $\widehat{\Gamma}$ l'ensemble des bons points pour la théorie de Pesin et le théorème d'Oseledets.

On a alors le résultat suivant:

\begin{Thm}{\label{codage}}

On considère une mesure $\mu$ comme ci-dessus avec $h_{\mu}(f)>0$.

Soit $\delta>0$. Pour tout $\gamma_1 >0$ et $\rho >0$, il existe un ensemble $\Lambda=\Lambda_{\delta, \rho, \gamma_1}$ avec $ \mu(\Lambda) \geq 1 - \delta$ tel que pour $y \in \Lambda$, on ait $\widehat{y} \in \widehat{\Gamma}$ avec $\pi(\widehat{y})=y$, $\eta(\widehat{y})>0$ et $n=n(\widehat{y}) > 1$ un entier avec le diagramme commutatif suivant: 

$$
\xymatrix{
\Sigma=\{1,\cdots,N\}^\mathbb{N} \ar[r]^{\sigma}  \ar[d]_{\sigma_v}  &\ar[d]^{\sigma_v} \Sigma= \{1,\cdots,N\}^\mathbb{N}   \\
G_v \ar[r]^{\Gamma} & G_v
}
$$

où $N= e^{h_{\mu}(f)n - \rho n}$, $\sigma$ est le décalage à gauche, $G_v$ est l'ensemble des graphes de fonctions holomorphes $(\phi(Y),Y)$ au-dessus de $\overline{B_{k-s}(0,\eta(\widehat{y}))}$ avec $\mathrm{Lip}( \phi) \leq\gamma_1$, $\sigma_v$ est un homéomorphisme de $\Sigma$ sur son image. L'application $\Gamma$ est définie et continue sur $\sigma_v(\Sigma)$ et elle est reliée à $f$ par:

$$\zeta \in \sigma_v(\Sigma) \Rightarrow C_{\gamma}^{-1}(\widehat{y}) \circ  \tau_y^{-1} \circ f^n \circ  \tau_y \circ C_{\gamma}(\widehat{y})(\zeta) \subset \Gamma(\zeta).$$

Soit $\beta= {\sigma_v}_* (\nu_0)$ où $\nu_0$ est la mesure de Bernoulli de $\Sigma$. Alors, $\beta$ est invariante par $\Gamma$ et a une entropie $h_{\beta}(\Gamma)=\log N$.

Enfin, soit $S_{\widehat{y}}$ le courant uniformément laminaire (voir \cite{BLS}) défini par $S_{\widehat{y}}= \int [\sigma_v(\omega)] d \nu_0(\omega)$, alors

$$\mbox{dim}_{\Hcal}(\mbox{Support de } S_{\widehat{y}}) \geq 2(k-s) + \frac{h_{\mu}(f)- \rho}{\chi_1 + 4 \gamma}.$$

\end{Thm}

Voici le plan de l'article. Dans un premier paragraphe, nous allons démontrer le théorème \ref{codage}. Nous verrons ensuite comment il implique le théorème \ref{general} et nous finirons par la démonstration du théorème \ref{bedford-diller}.

\section{\bf{Démonstration du théorème \ref{codage}}}

Soit $\delta > 0$. On fixe $\gamma>0$ petit devant les exposants de Lyapounov et on applique le théorème d'Oseledets avec ce $\gamma$.

Comme les exposants de Lyapounov ne dépendent pas du choix des cartes pour $X$, on peut supposer que $\mu$ ne charge pas le bord des cartes. En particulier, pour $\epsilon_1>0$ suffisamment petit, la masse pour $\mu$ d'un $\epsilon_1$-voisinage des bords des cartes qui recouvrent $X$ (pour une métrique fixée sur $X$) est plus petite que $\frac{\delta}{16}$.

On note $V_{\epsilon_1}$ ce voisinage et on a $\widehat{\mu}(\pi^{-1}(V_{\epsilon_1}^c))=\mu(V_{\epsilon_1}^c)\geq 1-\frac{\delta}{16}$.

Par le théorème de Lusin, on peut trouver un ensemble $\Lambda_{\delta}'$ compact, inclus dans $\pi^{-1}(V_{\epsilon_1}^c)$ avec $\widehat{\mu}(\Lambda_{\delta}')\geq 1- \frac{\delta}{8 }$ et $\widehat{x}\mapsto C_{\gamma}^{\pm 1}(\widehat{x})$, $\widehat{x}\mapsto r(\widehat{x})$ continues sur $\Lambda_{\delta}'$. 

Soit $\gamma_1 > 0$ et $\rho > 0$. On continue de reprendre les notations de \cite{DetNgu}.

Dans un premier temps, nous allons construire l'ensemble $\Lambda$ du théorème, puis nous montrerons ses propriétés.

\subsection{\bf{Construction de $\Lambda$}}

Soit $0< \gamma_0 < \gamma_1$ fixé suffisamment petit pour que $\gamma_0 < e^{\frac{\gamma}{4}}-1$, $\gamma_0 < 1-e^{- \frac{\gamma}{2}}$ et $\gamma_0e^{\chi_1+\gamma}+e^{-4\gamma} < e^{-\frac{7}{2}\gamma}$. 

$\Lambda_{\delta}'$ est compact, il existe donc $r_0>0$ tel que : 

$$\forall \widehat{x}\in \Lambda_{\delta}' ,\ r(\widehat{x})\geq r_0 \mbox{  et  } r_0\leq \| C_{\gamma}(\widehat{x})^{\pm 1}\| \leq \frac{1}{r_0}.$$

On fixe $0 < r < 1$ petit de sorte que $0<r<\rho$ et $\frac{h_{\mu}(f)}{1+r}-3r>h_{\mu}(f)-\rho$. 

Par le théorème de Brin-Katok, si on note

$$B_{m} (x, \epsilon)= \{ y \in X \mbox{  ,  } dist(f^{p}(x),f^{p}(y)) < \epsilon \mbox{  pour  } p=0, \cdots m-1 \},$$ 

on a

$$h_{\mu}(f)=\lim_{\epsilon\rightarrow 0} \liminf_{m \rightarrow + \infty} -\frac{1}{m} \log \mu (B_m(x,\epsilon)) \mbox{    } \mbox{  pour  }\mu-\mbox{presque tout  }x.$$

On en déduit que si on pose: 

$$\Gamma_{\epsilon,m_0}=\{x, \forall m\geq m_0\ \mu (B_m(x,\epsilon))\leq e^{-h_{\mu}(f)m+rm}\},$$

on a $\mu (\Gamma_{\epsilon,m_0})\geq 1-\frac{\delta}{4}$ pour $\epsilon$ assez petit puis $m_0$ suffisamment grand.

Dans la suite on considère le $h$ et le $\eta$ du Closing Lemma ou de sa preuve dans \cite{DetNgu}. 

Les $\{B(\widehat{y}, \eta/4),\widehat{y}\in \Lambda_{\delta}' \}$ forment un recouvrement de $\Lambda_{\delta}'$ qui est compact: on peut donc trouver un sous-recouvrement fini $\cup_{i=1}^{t} B(\widehat{y_{i}}, \eta/4)$. 

Soit $\alpha > 0$ tel que $\alpha t < \frac{\delta}{2}$.

On considère $\xi$ une partition finie de l'extension naturelle $\widehat{X}$, plus fine que $(\Lambda_{\delta}',\widehat{X}\setminus \Lambda_{\delta}')$ et telle que le diamètre des éléments de $\xi$ soit plus petit que $\eta/2$. 

On note: 

\begin{equation*}
\begin{split}
\Lambda_{\delta,m_0'} & =\{ \widehat{x}\in \Lambda_{\delta}', \forall m\geq m_0' \mbox{  on a  } \widehat{f}^q(\widehat{x})\in \xi(\widehat{x}) \mbox{  pour un  } q\in [m,m(1+r)] \}. \\
\end{split}
\end{equation*}

Comme dans \cite{KH} et \cite{DetNgu}, pour $m_0'$ assez grand, on a $\widehat{\mu}(\Lambda_{\delta,m_0'})\geq 1-\delta/4$.

Ainsi on a $\mu(\pi(\Lambda_{\delta,m_0'})\cap \Gamma_{\epsilon,m_0})\geq 1-\delta/2$.

Soit 

$$I_1=\{i \in \{1, \cdots ,t\} \mbox{  ,  } \mu(\pi(\Lambda_{\delta,m_0'} \cap B(\widehat{y_{i}}, \eta/4)) \cap \Gamma_{\epsilon,m_0}) < \alpha \}$$

et $I_2=I \setminus I_1$.

On a 

$$\mu \left( \cup_{i \in I_1} \pi(\Lambda_{\delta,m_0'} \cap B(\widehat{y_{i}}, \eta/4)) \cap \Gamma_{\epsilon,m_0} \right) < t \alpha < \frac{\delta}{2}.$$

Enfin, on notera

$$\Lambda = \Lambda_{\delta, \rho, \gamma_1} = \cup_{i \in I_2} \pi(\Lambda_{\delta,m_0'} \cap B(\widehat{y_{i}}, \eta/4)) \cap \Gamma_{\epsilon,m_0}.$$

C'est le $\Lambda$ que l'on voulait construire et on a $\mu(\Lambda) \geq 1 -\delta/2 - \delta/2=1- \delta$.

Montrons maintenant que les points de cet ensemble vérifient les propriétés du théorème.

\subsection{\bf{Propriétés de $\Lambda$}}

Soit $y \in \Lambda$. Par définition de cet ensemble, il existe $\widehat{y} \in \Lambda_{\delta,m_0'} \cap B(\widehat{y_{i}}, \eta/4)$ avec $\pi(\widehat{y})=y$ et $i \in I_2$.

Comme $B(\widehat{y_{i}}, \eta/4) \subset B(\widehat{y}, \eta/2)$, on a 

$$\mu ( \pi(\Lambda_{\delta,m_0'} \cap B(\widehat{y}, \eta/2)) \cap \Gamma_{\epsilon,m_0} ) \geq \alpha$$

par définition de $I_2$.

Soit $m\geq \max(m_0,m_0')$ de sorte que $e^{-rm} < \alpha$. On prend $E_m$ un ensemble $(m,\epsilon)$-séparé de cardinal maximal dans $\pi(\Lambda_{\delta,m_0'} \cap  B(\widehat{y}, \eta/2))\cap \Gamma_{\epsilon,m_0}$. 

Si $x\in E_m$, on a $\mu(B_m(x,\epsilon))\leq e^{-h_{\mu}(f)m+rm}$. Le cardinal de $E_m$ est donc supérieur à $\alpha e^{h_{\mu}(f)m-rm}$.

Pour $q=m,\cdots,[(1+r)m]$, on pose $V_q=\{x\in E_m, \widehat{f}^q(\widehat{x})\in \xi(\widehat{x})\}$ où $\widehat{x}$ est un point de $\Lambda_{\delta,m_0'} \cap  B(\widehat{y}, \eta/2)$ qui se projette sur $x$. 

Soit $n$ qui maximise le cardinal de $V_q$. On a:

$$\mathrm{card}\ V_n \geq \frac{\alpha e^{h_{\mu}(f)m-rm}}{rm+1}\geq \alpha e^{h_{\mu}(f)m-2rm} \geq e^{h_{\mu}(f)m-3rm} \geq e^{h_{\mu}(f)\frac{n}{1+r}-3rn}  \geq e^{h_{\mu}(f)n-\rho n}.$$

Notons $x_1 , \cdots , x_N$ les points de $V_n$. Quitte à réduire leur nombre, on pourra supposer que $N= e^{h_{\mu}(f)n-\rho n}$.

Par construction $x_1,\cdots,x_N$ sont dans $\pi(\Lambda_{\delta,m_0'} \cap  B(\widehat{y}, \eta/2))$ et nous appellerons $\widehat{x_1},\cdots,\widehat{x_N}$ les points de $\Lambda_{\delta,m_0'} \cap  B(\widehat{y}, \eta/2)$ qui se projettent sur eux. 

Pour $i=1, \cdots , N$, on a $\widehat{f}^n( \widehat{x_i}) \in \xi(\widehat{x_i})$ ce qui implique que $d(\widehat{f}^n( \widehat{x_i}),\widehat{x_i}) < \eta/2$. En particulier, $\widehat{x_i}$ et $\widehat{f}^n( \widehat{x_i})$ sont dans $B(\widehat{y}, \eta)$.

C'est exactement la même situation qu'à la fin du paragraphe 3.1 de \cite{DetNgu}. Nous pouvons donc utiliser ce que nous avons fait dans \cite{DetNgu} en particulier la construction des graphes verticaux qui sont codés par $\Sigma=\{1, \cdots , N \}^{\Nn}$. En voici rapidement leurs constructions.

On se place dans le repère $C_{\gamma}^{-1}(\widehat{y})E_u(\widehat{y})\oplus C_{\gamma}^{-1}(\widehat{y})E_s(\widehat{y})$.

Soit $\mathcal{A}^0$ l'ensemble constitué de l'union des graphes d'applications holomorphes $(\psi(Y),Y)$ au-dessus de $C_{\gamma}^{-1}(\widehat{y})E_s(\widehat{y})$ avec $ \| Y\|\leq hr(\widehat{y})e^{\frac{\gamma}{2}}$, $\psi(Y)=\mathrm{constante}$,  $\| \psi(0)\|\leq hr(\widehat{y})e^{-\frac{\gamma}{2}}$ et $\mathrm{Lip}( \psi)=0 \leq \gamma_0e^{-\frac{\gamma}{2}}$.

On prend l'image de ces graphes par $C_i=C_{\gamma}^{-1}(\widehat{f}^n(\widehat{x_i}))\tau_{f^n(x_i)}^{-1}\tau_y C_{\gamma}(\widehat{y})$ et on obtient des graphes $(\phi(Y),Y)$ au-dessus d'une partie de $C_{\gamma}^{-1}(\widehat{f}^n(\widehat{x_i}))E_s(\widehat{f}^n(\widehat{x_i}))$ puis on tire en arrière ces graphes $(\phi(Y),Y)$ par $g_{\widehat{f}^{n-1}(\widehat{x_i})},\cdots,g_{\widehat{x_i}}$. On obtient à  la fin des graphes $(\phi_0(Y),Y)$ au-dessus d'une partie de $C_{\gamma}^{-1}(\widehat{x_i})E_s(\widehat{x_i})$ pour au-moins $\| Y\|\leq e^{\gamma}hr(\widehat{x_i})$, $\| \phi_0(0)\|\leq e^{-\gamma}hr(\widehat{x_i})$ et $\mathrm{Lip}(\phi_0)\leq e^{-\gamma}\gamma_0$. 

Maintenant, on remet ces graphes dans le repère initial $C_{\gamma}^{-1}(\widehat{y})E_u(\widehat{y})\oplus C_{\gamma}^{-1}(\widehat{y})E_s(\widehat{y})$ c'est-à-dire que l'on prend leur image par $C_{\gamma}^{-1}(\widehat{y})\tau_y^{-1}\tau_{x_i} C_{\gamma}(\widehat{x_i})$. On obtient des graphes $(\psi_0(Y),Y)$ pour au moins $\| Y\|\leq e^{\frac{\gamma}{2}}hr(\widehat{y})$ avec $\| \psi_0(0)\| \leq e^{-\frac{3\gamma}{4}}hr(\widehat{y})$ et $\mathrm{Lip} (\psi_0) \leq e^{-\frac{\gamma}{2}}\gamma_0$. Ce sont les mêmes conditions qu'au départ (c'est-à-dire que l'on pourra recommencer le procédé avec un $\widehat{x_j}$ pas nécessairement égal à $\widehat{x_i}$). 

Notons $\mathcal{A}_1^0,\cdots,\mathcal{A}_N^0$ les ensembles constitués de l'union de ces graphes obtenus pour $i=1,\cdots,N$. On a $\mathcal{A}_i^0\subset \mathcal{A}^0$. On notera aussi $\mathcal{A}_i$ l'image de $\mathcal{A}_i^0$ par $\tau_y\circ C_{\gamma}(\widehat{y})$. 

A la $(l+1)$-ème génération on a $N^{l+1}$ ensembles $\mathcal{A}_{w_0\cdots w_{-l}}^0$ avec $w_i=1,\cdots,N$. 

$\mathcal{A}_{w_0\cdots w_{-l}}^0$ est constitué de graphes $(\psi(Y),Y)$ avec $\| Y\|\leq e^{\frac{\gamma}{2}}hr(\widehat{y})$, $\| \psi(0)\|\leq e^{-\frac{3\gamma}{4}}h r(\widehat{y})$ et $\mathrm{Lip}( \psi) \leq \gamma_0 e^{-\frac{\gamma}{2}}$. 

Ces graphes proviennent des graphes de $\mathcal{A}_{w_{-1}\cdots w_{-l}}^0$ via le procédé précédent qui utilise $g_{\widehat{x_{w_0}}}$ , $\cdots$ ,$g_{\widehat{f}^{n-1}(\widehat{x_{w_0}})}$. 

On a vu dans \cite{DetNgu} que $\mathcal{A}_{w_{0}\cdots w_{-l}}^0\subset \mathcal{A}_{w_{0}\cdots w_{-l+1}}^0$ pour tout $w_0,\cdots,w_{-l}\in \{1,\cdots,N\}$ et $l\geq 0$ et que les $N^{l+1}$ ensembles $\overline{\mathcal{A}_{w_0\cdots w_{-l}}^0}$ sont disjoints. 

Les graphes qui constituent $\mathcal{A}_{w_{0}\cdots w_{-l}}^0$ sont dans $G_v$ qui est l'ensemble des graphes de fonctions holomorphes $(\phi(Y),Y)$ au-dessus de $\overline{B_{k-s}(0,\eta(\widehat{y}))}$ avec $\mathrm{Lip}( \phi) \leq \gamma_1$ et $\eta(\widehat{y})=e^{\gamma/2}hr(\widehat{y})$.

Sur $G_v$ on met la métrique: 

$$\mathrm{d}(\mathrm{graphe} \mbox{    } \phi,\mathrm{graphe} \mbox{    } \psi)=\max_{\overline{B_{k-s}(0,\eta(\widehat{y}))}}\| \phi-\psi\|.$$

On peut maintenant définir une application de $\{1,\cdots,N\}^{\mathbb{N}}$ dans $G_v$ de la façon suivante: soit $w=(w_0,w_{-1},\cdots,w_{-l},\cdots)\in \{1,\cdots,N\}^{\mathbb{N}}$. On prend $A_l$ un graphe qui constitue $\mathcal{A}_{w_0 w_{-1}\cdots w_{-l}}^0$ pour $l\geq 0$. 

$A_l$ est dans $G_v$ et par le théorème d'Ascoli, il existe une sous-suite $(A_{l_j})$ qui converge vers $A \in G_v$ (une limite uniforme de fonctions holomorphes est holomorphe). On a montré dans \cite{DetNgu} que $(A_l)_l$ converge aussi vers $A$ (et que cet ensemble est indépendant du choix de $A_l$ comme graphe qui constitue $\mathcal{A}_{w_0 w_{-1}\cdots w_{-l}}^0$). On a ainsi défini une application: 

\begin{center}
\begin{tabular}{cccc}
$\sigma_v:$ & $\{1,\cdots,N\}^{\mathbb{N}}$ & $\longrightarrow$ & $G_v$\\
             & $w$                         & $\longrightarrow$ & $A$ 
\end{tabular}
\end{center}

On a vu dans \cite{DetNgu}, que cette application $\sigma_v$ est injective et que l'on a même $\sigma_v(w)\cap \sigma_v(w')=\varnothing$ si $w\neq w'$. 

Si  $\sigma$ désigne le décalage à gauche, on peut donc définir l'application $\Gamma$ par

$$\Gamma(\zeta):=\Gamma(\sigma_v(\omega))=\sigma_v(\sigma(\omega))$$

pour $\zeta \in \sigma_v(\Sigma)$.

On a donc obtenu le diagramme commutatif

$$
\xymatrix{
\Sigma=\{1,\cdots,N\}^\mathbb{N} \ar[r]^{\sigma}  \ar[d]_{\sigma_v}  &\ar[d]^{\sigma_v} \Sigma= \{1,\cdots,N\}^\mathbb{N}   \\
G_v \ar[r]^{\Gamma} & G_v
}
$$

Il faut maintenant montrer les propriétés énoncées dans le théorème. Sur $\Sigma$, on met la métrique $d(\omega, \omega')= \sum_{i=0}^{+ \infty} \frac{| \omega_i - \omega_i'|}{2^{i}}$.

Tout d'abord $\sigma_v$ est continue. En effet, soit $\beta >0$. 

Pour $p$ assez grand, on a $2he^{-2 \gamma n (p+1) + 2 \gamma (p+1)} \leq \beta$.

Si $d(\omega, \omega') < \frac{1}{2^{p+1}}$, avec $\omega, \omega' \in \Sigma$, on a $\omega_0 = \omega_0' , \cdots , \omega_{-p} = \omega_{-p}'$.

Par définition, $\sigma_v(\omega)= \lim A_l^1$ avec $A_l^1$ un graphe qui constitue $\mathcal{A}_{w_{0}\cdots w_{-l}}^0$ et $\sigma_v(\omega')= \lim A_l^2$ avec $A_l^2$ un graphe de $\mathcal{A}_{w_{0}'\cdots w_{-l}'}^0$.

Mais pour $l \geq p$, on a $\mathcal{A}_{w_{0}\cdots w_{-l}}^0 \subset \mathcal{A}_{w_{0}\cdots w_{-p}}^0$ et $\mathcal{A}_{w_{0}'\cdots w_{-l}'}^0 \subset \mathcal{A}_{w_{0}'\cdots w_{-p}'}^0 =\mathcal{A}_{w_{0}\cdots w_{-p}}^0 $ d'où 

$$d(A_l^1,A_l^2) \leq 2he^{-2 \gamma n (p+1) + 2 \gamma (p+1)} \leq \beta$$

par le lemme 3.6 de \cite{DetNgu} car chaque point de $A_l^1$ et $A_l^2$ appartient à un graphe qui constitue $\mathcal{A}_{w_{0}\cdots w_{-p}}^0$. 

Par passage à la limite, on a $d(\sigma_v(\omega),\sigma_v(\omega')) \leq \beta$ ce qui montre que $\sigma_v$ est continue.

L'application $\sigma_v$ est une bijection continue de $\Sigma$ sur $\sigma_v(\Sigma)$ c'est donc un homéomorphisme car $\Sigma$ est un espace métrique compact.

Par ailleurs, sur $\sigma_v(\Sigma)$ on a $\Gamma = \sigma_v \circ \sigma \circ \sigma_v^{-1}$ d'où la continuité de $\Gamma$ aussi.

Etablissons maintenant le lien entre $\Gamma$ et $f$.

Soit $\omega=(\omega_0 , \cdots , \omega_{-l} , \cdots) \in \Sigma$. On a $\sigma_v(\omega)=A$ où $A= \lim A_l^0$ avec $A_l^0$ un graphe qui constitue $\mathcal{A}_{w_{0}\cdots w_{-l}}^0$.

Pour $l\geq 1$, on a par construction (voir le début du paragraphe 3.3 de \cite{DetNgu}) $$C_{\gamma}^{-1}(\widehat{y})\tau_y^{-1}\tau_{f^n(x_{w_0})}C_{\gamma}(\widehat{f}^n(\widehat{x_{w_0}}))g_{\widehat{f}^{n-1}(\widehat{x_{w_0}})}\circ \cdots \circ g_{\widehat{x_{w_0}}}\circ C_{\gamma}^{-1}(\widehat{x_{w_0}})\tau_{x_{w_0}}^{-1}\tau_y C_{\gamma}(\widehat{y})(A_l^0)$$

qui est inclus dans un graphe $A_l^1$ de $\mathcal{A}_{w_{-1}\cdots w_{-l}}^0$.

Comme $(A_l^1)$ converge vers $\sigma_v(( \omega_{-1}, \cdots , \omega_{-l}, \cdots ))$, par passage à la limite on a

$$C_{\gamma}^{-1}(\widehat{y})\tau_y^{-1}\tau_{f^n(x_{w_0})}C_{\gamma}(\widehat{f}^n(\widehat{x_{w_0}}))g_{\widehat{f}^{n-1}(\widehat{x_{w_0}})}\circ \cdots \circ g_{\widehat{x_{w_0}}}\circ C_{\gamma}^{-1}(\widehat{x_{w_0}})\tau_{x_{w_0}}^{-1}\tau_y C_{\gamma}(\widehat{y})(\sigma_v(\omega))$$

qui est inclus dans $\sigma_v(\sigma(\omega))$, c'est-à-dire

$$C_{\gamma}^{-1}(\widehat{y}) \circ \tau_y^{-1} \circ f^n \circ \tau_y \circ C_{\gamma}(\widehat{y})(\sigma_v(\omega)) \subset \sigma_v(\sigma(\omega))=\Gamma(\sigma_v(\omega)).$$

C'est ce que l'on voulait montrer.

Soit $\beta= {\sigma_v}_* \nu_0$ où $\nu_0$ est la mesure de Bernoulli de $\Sigma$. On a 

$$\Gamma_* \beta= \Gamma_* {\sigma_v}_* \nu_0 = {\sigma_v}_* \sigma_* \nu_0 = {\sigma_v}_*  \nu_0 = \beta$$

c'est-à-dire que la mesure $\beta$ est invariante par $\Gamma$.

De plus, comme $\sigma_v$ est un homéomorphisme, on a $h_{\beta}(\Gamma)=h_{\nu_0}(\sigma)= \log N$.

Considérons enfin le courant uniformément laminaire (voir \cite{BLS}) $S_{\widehat{y}}= \int [\sigma_v(\omega)] d \nu_0(\omega)$: le but maintenant est de minorer la dimension de Hausdorff de son support.

Si $\omega=(\omega_0, \cdots , \omega_{-l} , \cdots)$, on a $\sigma_v(\omega) \subset \overline{\mathcal{A}_{w_{0}\cdots w_{-l}}^0}$ pour tout $l \geq 0$. En particulier,

$$\mbox{Support de } S_{\widehat{y}} \subset  \bigcup_{\omega_0, \cdots , \omega_{-l}=1, \cdots ,N} \overline{\mathcal{A}_{w_{0}\cdots w_{-l}}^0}.$$

Soit $\alpha_0$ la distance minimale entre deux $\overline{\tau_y \circ C_{\gamma}(\widehat{y})(\mathcal{A}_i^0)}$ ($i=1, \cdots , N$). Comme les $\mathcal{A}_i^0$ sont disjoints (voir le lemme 3.3 de \cite{DetNgu}), on a $\alpha_0 > 0$.

On considère $x$ dans le support de $S_{\widehat{y}}$ et on va majorer  $S_{\widehat{y}} \wedge \Omega_0^s (B(x, e^{- \chi_1 nl - 4 \gamma nl}))$ pour $l$ grand (ici $\Omega_0$ est la $(1,1)$-forme standard de $\Cc^k$).

Si $\zeta$ est une lame de $S_{\widehat{y}}$, par construction c'est un graphe $(\Phi(Y),Y)$ pour $Y \in B_{k-s}(0, \eta(\widehat{y}))$ avec $\mathrm{Lip}( \Phi) \leq \gamma_1$. Par la formule de la coaire (voir \cite{Fe} p.258), on a donc

$$\zeta \wedge \Omega_0^s (B(x, e^{- \chi_1 nl - 4 \gamma nl})) \leq e^{2(k-s)(- \chi_1 nl - 4 \gamma nl)})(1+ \gamma_1)^{2(k-s)}$$

et on obtient donc une majoration de $S_{\widehat{y}} \wedge \Omega_0^s (B(x, e^{- \chi_1 nl - 4 \gamma nl}))$ par

$$ e^{2(k-s)(- \chi_1 nl - 4 \gamma nl)})(1+ \gamma_1)^{2(k-s)} \nu_0( \{ \omega , \sigma_v(\omega) \cap B(x, e^{- \chi_1 nl - 4 \gamma nl}) \neq \emptyset \} ) .$$

Maintenant, on a 

\begin{Lem}{\label{lemme1}}

Pour $l$ assez grand par rapport à $\alpha_0$, $B(x, e^{- \chi_1 nl - 4 \gamma nl})$ ne rencontre qu'un seul $\overline{\mathcal{A}_{w_{0}\cdots w_{-l}}^0}$. Ici le $l$ est uniforme en $x$ sur le support de $S_{\widehat{y}}$.

\end{Lem}

Admettons ce lemme pour le moment et terminons la minoration de la dimension du support de $S_{\widehat{y}}$.

Soit $\overline{\mathcal{A}_{w_{0}' \cdots w_{-l}'}^0}$ le seul des ensembles ci-dessus qui rencontre $B(x, e^{- \chi_1 nl - 4 \gamma nl})$. On a donc

\begin{equation*}
\begin{split}
S_{\widehat{y}} \wedge \Omega_0^s (B(x, e^{- \chi_1 nl - 4 \gamma nl})) & \leq  
e^{2(k-s)(- \chi_1 nl - 4 \gamma nl)})(1+ \gamma_1)^{2(k-s)} \nu_0( \{ \omega \mbox{ ,  } \sigma_v(\omega) \subset \overline{\mathcal{A}_{w_{0}' \cdots w_{-l}'}^0} \} ) \\
& \leq e^{2(k-s)(- \chi_1 nl - 4 \gamma nl)})(1+ \gamma_1)^{2(k-s)} \nu_0( \{ \omega \mbox{ ,  } \omega_0 = \omega_0' \cdots \omega_{-l} = \omega_{-l}' \}) \\
& \leq e^{2(k-s)(- \chi_1 nl - 4 \gamma nl)})(1+ \gamma_1)^{2(k-s)} \left( \frac{1}{N} \right)^{l+1}\\
& =   e^{2(k-s)(- \chi_1 nl - 4 \gamma nl)})(1+ \gamma_1)^{2(k-s)}  e^{(-h_{\mu}(f)n + \rho n)(l+1)}.
\end{split}
\end{equation*}

Ainsi, en utilisant la remarque de la proposition 2.1 de \cite{Y}, 

$$\underline{d}_{S_{\widehat{y}} \wedge \Omega_0^s}(x)= \liminf_{l \rightarrow + \infty} \frac{\log S_{\widehat{y}} \wedge \Omega_0^s (B(x, e^{- \chi_1 nl - 4 \gamma nl}))}{\log e^{- \chi_1 nl - 4 \gamma nl}} \geq 2(k-s) + \frac{h_{\mu}(f) - \rho }{\chi_1 + 4 \gamma}. $$ 

Comme cette inégalité est vraie pour tout $x$ dans le support de $S_{\widehat{y}}$, on a par la proposition 2.1 de \cite{Y} que

$$\mbox{dim}_{\Hcal}(\mbox{Support de } S_{\widehat{y}}) \geq 2(k-s) + \frac{h_{\mu}(f)- \rho}{\chi_1 + 4 \gamma}.$$

Remarquons qu'à ce stade le courant $S_{\widehat{y}}$ dépend de $n$ et donc de $\gamma$ et $\rho$. 

Signalons aussi que l'inégalité précédente peut s'exprimer en terme de dimension de courant (voir \cite{DV}).

Passons maintenant à la preuve du lemme précédent.

\bigskip

{\bf Démonstration du lemme \ref{lemme1}}

On considère $l$ suffisamment grand pour que $e^{- \gamma n l} \leq h r_0 e^{- \gamma n}$. 

Le point $x$ est dans le support de $S_{\widehat{y}}$ et il est donc inclus dans un $\overline{\mathcal{A}_{w_{0}\cdots w_{-l}}^0}$. Pour $(w_{0}', \cdots , w_{-l}') \neq (w_{0}, \cdots , w_{-l})$, montrons donc que $\overline{\mathcal{A}_{w_{0}' \cdots w_{-l}'}^0}$ ne rencontre pas $B(x, e^{- \chi_1 nl - 4 \gamma nl})$.

Si ce n'est pas le cas, soit $z$ dans $\overline{\mathcal{A}_{w_{0}' \cdots w_{-l}'}^0} \cap B(x, e^{- \chi_1 nl - 4 \gamma nl})$.

Par le lemme 3.9 de \cite{DetNgu}, les points $\tau_y C_{\gamma}(\widehat{y})(x)$ et $\tau_y C_{\gamma}(\widehat{y})(z)$ sont $(l+1, \alpha_0)$-séparés pour $f^n$. 

Comme $d(\widehat{x_{\omega_0}}, \widehat{y}) < \eta$, par \cite{DetNgu}, on a pour $t$ dans le segment $[xz]$

$$\| D(C_{\gamma}^{-1}(\widehat{x_{\omega_0}}) \tau_{x_{\omega_0}}^{-1} \tau_y C_{\gamma}(\widehat{y}))(t)\|= \| C_{\gamma}^{-1}(\widehat{x_{\omega_0}}) \circ C_{\gamma}(\widehat{y}) \| \leq e^{\gamma}.$$

La distance entre les points $C_{\gamma}^{-1}(\widehat{x_{\omega_0}}) \tau_{x_{\omega_0}}^{-1} \tau_y C_{\gamma}(\widehat{y})(x)$ et $C_{\gamma}^{-1}(\widehat{x_{\omega_0}}) \tau_{x_{\omega_0}}^{-1} \tau_y C_{\gamma}(\widehat{y})(z)$ est donc majorée par $e^{- \chi_1 nl - 4 \gamma nl + \gamma}$.

Maintenant les points 

$g_{\widehat{f}^{p-1}(\widehat{x_{w_0}})}\circ \cdots \circ g_{\widehat{x_{w_0}}}\circ C_{\gamma}^{-1}(\widehat{x_{w_0}})\tau_{x_{w_0}}^{-1}\tau_y C_{\gamma}(\widehat{y})(x)$ sont dans les boîtes $\overline{B_{s}(0,e^{\gamma} hr(\widehat{f}^p(\widehat{x_{w_0}})))}\times \overline{B_{k-s}(0,e^{\gamma} hr(\widehat{f}^p(\widehat{x_{w_0}})))} \subset B_k(0, 2 h r(\widehat{f}^p(\widehat{x_{w_0}}))) $ pour $p=0,\cdots,n-1$ par construction des graphes qui composent $\mathcal{A}_{w_{0}\cdots w_{-l}}^0$ (voir le début du paragraphe 3.3 de \cite{DetNgu}).

Sur $B(0, 3 h r(\widehat{f}^p(\widehat{x_{w_0}}))) $, on a par la proposition 1.1 de \cite{DetNgu}, pour $h$ assez petit,

$$\| Dg_{\widehat{f}^{p}(\widehat{x_{w_0}})}(w) \| \leq e^{\chi_1 + \gamma} + \frac{ \|w\|}{r(\widehat{f}^p(\widehat{x_{w_0}}))} \leq e^{\chi_1 + 2 \gamma}.$$

Ainsi, par récurrence sur $p=0, \cdots , n-1$ et par l'inégalité des accroissements finis, l'image du segment $[xz]$ par $g_{\widehat{f}^{n-1}(\widehat{x_{w_0}})}\circ \cdots \circ g_{\widehat{x_{w_0}}}\circ C_{\gamma}^{-1}(\widehat{x_{w_0}})\tau_{x_{w_0}}^{-1}\tau_y C_{\gamma}(\widehat{y})$ a une longueur plus petite que $e^{- \chi_1 nl - 4 \gamma nl + \chi_1 n + 2 \gamma n+  \gamma}$. Ici, on utilise que $h r(\widehat{f}^p(\widehat{x_{w_0}})) \geq h r_0 e^{- \gamma n} \geq e^{- \gamma n l}$ pour $p=0, \cdots , n-1$.

Maintenant, si on prend l'image par $C_{\gamma}^{-1}(\widehat{y}) \tau_{y}^{-1} \tau_{f^{n}(x_{w_0})} C_{\gamma}(\widehat{f}^n(\widehat{x_{w_0}}))$, on obtient que l'image de $[xz]$ par $C_{\gamma}^{-1}(\widehat{y}) \tau_{y}^{-1} \tau_{f^{n}(x_{w_0})} C_{\gamma}(\widehat{f}^n(\widehat{x_{w_0}})) g_{\widehat{f}^{n-1}(\widehat{x_{w_0}})}\circ \cdots \circ g_{\widehat{x_{w_0}}}\circ C_{\gamma}^{-1}(\widehat{x_{w_0}})\tau_{x_{w_0}}^{-1}\tau_y C_{\gamma}(\widehat{y})$ a une longueur plus petite que $e^{- \chi_1 nl - 4 \gamma nl + \chi_1 n + 2 \gamma n+  2 \gamma}$. 

On recommence ce que l'on vient de faire avec $\widehat{x_{w_{-1}}}$ à la place de $\widehat{x_{w_0}}$ et ainsi de suite. On obtient donc que l'image de $[xz]$ par 

\begin{equation*}
\begin{split}
&C_{\gamma}^{-1}(\widehat{y}) \tau_{y}^{-1} \tau_{f^{n}(x_{w_{-p+1}})} C_{\gamma}(\widehat{f}^n(\widehat{x_{w_{-p+1}}})) g_{\widehat{f}^{n-1}(\widehat{x_{w_{-p+1}}})}\circ \cdots \circ g_{\widehat{x_{w_{-p+1}}}}\circ C_{\gamma}^{-1}(\widehat{x_{w_{-p+1}}})\tau_{x_{w_{-p+1}}}^{-1}
 \cdots\\
& \cdots \tau_{f^{n}(x_{w_0})} C_{\gamma}(\widehat{f}^n(\widehat{x_{w_0}})) g_{\widehat{f}^{n-1}(\widehat{x_{w_0}})}\circ \cdots \circ g_{\widehat{x_{w_0}}}\circ C_{\gamma}^{-1}(\widehat{x_{w_0}})\tau_{x_{w_0}}^{-1}\tau_y C_{\gamma}(\widehat{y}) \\
\end{split}
\end{equation*}

a une longueur plus petite que $e^{- \chi_1 nl - 4 \gamma nl + \chi_1 np + 2 \gamma np+  2p \gamma}$ pour $p=0, \cdots , l$. Cette quantité est très petite devant $\alpha_0$ si $l$ est assez grand.

Mais la composée de l'application ci-dessus par $\tau_y C_{\gamma}(\widehat{y})$ au but est égale à $f^{np} \tau_y C_{\gamma}(\widehat{y})$ pour $p=0, \cdots , l$ et cela contredit le fait que les points  $\tau_y C_{\gamma}(\widehat{y})(x)$ et $\tau_y C_{\gamma}(\widehat{y})(z)$ sont $(l+1, \alpha_0)$-séparés pour $f^n$.

Cela termine la démonstration du lemme \ref{lemme1}.

\section{\bf{Démonstration du théorème \ref{general}}}

On applique le théorème \ref{codage} avec $\delta= 1/2$, $0 < \rho < \alpha_2$ et $\gamma_1=1/2$. 

Dans la suite, on fixe un $y \in \Lambda$ et on va montrer que le support de ${\tau_y}_* C_{\gamma}(\widehat{y})_* S_{\widehat{y}}$ (où $S_{\widehat{y}}$ est le courant uniformément laminaire obtenu dans le théorème \ref{codage}) est inclus dans le support de $T^+$. On obtiendra donc le résultat en faisant tendre $\rho$, puis $\gamma$ vers $0$ (car $T^+$ est indépendant de $\rho$ et $\gamma$).

Soit $x \notin \mbox{Supp}(T^+)$. On choisit $U \Subset \widetilde{U}$ des voisinages ouverts de $x$ tels que $T^+ \wedge \omega^{k-s}( \widetilde{U})=0$.

On prend $0 \leq \chi \leq 1$ une fonction lisse qui vaut $1$ sur $U$ et $0$ en dehors de $\widetilde{U}$. On a donc $\int \chi T^+ \wedge \omega^{k-s}=0$. On va montrer que $\int \chi {\tau_y}_* C_{\gamma}(\widehat{y})_* S_{\widehat{y}} \wedge \omega^{k-s}=0$. Cela impliquera que ${\tau_y}_* C_{\gamma}(\widehat{y})_* S_{\widehat{y}}(U)=0$ et donc que $x \notin \mbox{Supp}({\tau_y}_* C_{\gamma}(\widehat{y})_* S_{\widehat{y}})$, ce qui est le but.

On reprend la construction de $S_{\widehat{y}}$ . On rappelle que l'on part de $\mathcal{A}^0$ l'ensemble constitué de l'union des graphes d'applications holomorphes $(\psi(Y),Y)$ au-dessus de $C_{\gamma}^{-1}(\widehat{y})E_s(\widehat{y})$ avec $ \| Y\|\leq hr(\widehat{y})e^{\frac{\gamma}{2}}$, $\psi(Y)=\mathrm{constante}$,  $\| \psi(0)\|\leq hr(\widehat{y})e^{-\frac{\gamma}{2}}$ et $\mathrm{Lip}( \psi)=0 \leq \gamma_0e^{-\frac{\gamma}{2}}$.

Soit $\tau$ la transversale $E_u(\widehat{y}) \times \{0\}$. On note $\Delta^0(t)$ l'unique graphe de $\mathcal{A}^0$ qui passe par $(t,0)$ avec $\|t\| \leq hr(\widehat{y})e^{-\frac{\gamma}{2}}$. On définit $\beta^s= \int_{\tau} [\Delta^0(t)] d \lambda(t)$ où $\lambda$ est la mesure de Lebesgue sur $\tau$ normalisée ($\lambda(\tau)=1$). C'est une forme lisse et il existe une constante $C_1 > 0$ telle que $\beta^s \leq C_1 C_{\gamma}(\widehat{y})^* \tau_y^* (\omega^s)$.

Pour $l \geq 0$, chaque $\Delta^0(t)$ donne un graphe qui constitue $\mathcal{A}^0_{ \omega_0 \cdots \omega_{-l}}$ pour $\omega_0 , \cdots ,\omega_{-l}=1, \cdots , N$ (via le processus de tirés en arrière). On notera $\Delta^0_{\omega_0 \cdots \omega_{-l}}(t)$ le graphe en question. On a alors

\begin{Lem}

$$\int \frac{1}{N^{l+1}} \sum_{\omega_0 , \cdots ,\omega_{-l}=1, \cdots , N} [\Delta^0_{\omega_0 \cdots \omega_{-l}}(t)] d \lambda(t) \xrightarrow[l \longrightarrow + \infty]{} S_{\widehat{y}}.$$

\end{Lem}

\begin{proof}

Pour $l \geq 1$, soit

$$\nu_l= \frac{1}{N^{l+1}} \sum_{\omega_0 , \cdots ,\omega_{-l}=1, \cdots , N} \delta_{ \{ \omega' \in \Sigma \mbox{  ,  } \omega_0'=\omega_0 , \cdots ,\omega_{-l}'= \omega_{-l} \} }.$$

La suite $(\nu_l)$ converge vers $\nu_0$ qui est la mesure de Bernoulli sur $\Sigma$.

Soit $\Psi$ maintenant une (s,s) forme lisse. Montrons d'abord que l'application $\Sigma=\{ 1, \cdots , N \}^{\Nn} \longrightarrow \Cc$ qui à $\omega$ associe $\int_{\sigma_v(\omega)} \Psi$ est continue.

Soit $\epsilon > 0$. Si $\omega$ et $\omega'$ sont dans $\Sigma$ avec $d(\omega , \omega') < \frac{1}{2^{p+1}}$, on a vu dans la preuve du théorème \ref{codage} que $d(\sigma_v(\omega), \sigma_v(\omega')) \leq 2h e^{- \gamma n(p+1) + 2 \gamma (p+1)}$.

Dans la construction des $\sigma_v(\omega)$, on peut prendre un peu de marge quitte à réduire le $\eta(\widehat{y})$, c'est-à-dire que l'on peut avoir l'inégalité ci-dessus sur $B(0, e^{3 \gamma /4} h r(\widehat{y}))$ au lieu de $B(0, e^{\gamma/2 } h r(\widehat{y}))= B(0, \eta(\widehat{y}))$ si on veut. En particulier, grâce aux inégalités de Cauchy, si $\sigma_v(\omega)$ est un graphe $(\Phi_1(Y),Y)$ et $\sigma_v(\omega')$ un graphe $(\Phi_2(Y),Y)$, on a

$$\| D \Phi_1(Y) - D \Phi_2(Y) \| \leq \frac{2h e^{- \gamma n(p+1) + 2 \gamma (p+1)}}{(e^{3 \gamma /4}- e^{\gamma/2 }) h r(\widehat{y})} \leq \frac{2 e^{- \gamma n(p+1) + 2 \gamma (p+1)}}{(e^{3 \gamma /4}- e^{\gamma/2 }) r_0} $$

pour tout $Y \in  B(0, e^{\gamma/2 } h r(\widehat{y}))$.

On en déduit donc que

\begin{equation*}
\begin{split}
\left| \int_{\sigma_v(\omega)} \Psi - \int_{\sigma_v(\omega')} \Psi \right| &= \left| \int_{B(0, e^{\gamma/2 } h r(\widehat{y}))} [ (\Phi_1,Id)^* \Psi - (\Phi_2,Id)^* \Psi] \right| \\
& \leq C( \Psi, \gamma_1) \eta(\widehat{y})^{2(k-s)} \frac{2 e^{- \gamma n(p+1) + 2 \gamma (p+1)}}{(e^{3 \gamma /4}- e^{\gamma/2 }) r_0} \leq \epsilon \\
\end{split}
\end{equation*}

pour $p$ assez grand (indépendamment de $\omega$ et $\omega'$). On a donc montré que $\omega \longrightarrow \int_{\sigma_v(\omega)} \Psi$ est uniformément continue. Maintenant, on a

\begin{equation*}
\begin{split}
&\left|\langle S_{\widehat{y}}, \Psi \rangle - \int \frac{1}{N^{l+1}} \sum_{\omega_0 , \cdots ,\omega_{-l}=1, \cdots , N} \int_{\Delta^0_{\omega_0 \cdots \omega_{-l}}(t)} \Psi d \lambda(t) \right| \\
&= \left| \int \int_{\sigma_v(\omega)} \Psi d \nu_0(\omega) - \int \frac{1}{N^{l+1}} \sum_{\omega_0 , \cdots ,\omega_{-l}=1, \cdots , N} \int_{\Delta^0_{\omega_0 \cdots \omega_{-l}}(t)} \Psi d \lambda(t) \right| \\
& \leq \left| \int \int_{\sigma_v(\omega)} \Psi d \nu_0(\omega) -  \int \int_{\sigma_v(\omega)} \Psi d \nu_l(\omega) \right| \\
&+ \left| \int \int_{\sigma_v(\omega)} \Psi d \nu_l(\omega) - \int \frac{1}{N^{l+1}} \sum_{\omega_0 , \cdots ,\omega_{-l}=1, \cdots , N} \int_{\Delta^0_{\omega_0 \cdots \omega_{-l}}(t)} \Psi d \lambda(t) \right|. \\
\end{split}
\end{equation*}

Le premier terme tend vers $0$ car $\nu_l \longrightarrow \nu_0$ et $\omega \longrightarrow \int_{\sigma_v(\omega)} \Psi$ est continue.

Pour le second, soit $B(\omega_0, \cdots , \omega_{-l})=\{ \omega' \in \Sigma \mbox{  ,  } \omega_0'= \omega_0, \cdots , \omega_{-l}'= \omega_{-l} \}$. On a

\begin{equation*}
\begin{split}
& \left| \int \int_{\sigma_v(\omega')} \Psi d \nu_l(\omega') - \int \frac{1}{N^{l+1}} \sum_{\omega_0 , \cdots ,\omega_{-l}=1, \cdots , N} \int_{\Delta^0_{\omega_0 \cdots \omega_{-l}}(t)} \Psi d \lambda(t) \right| \\
&= \left| \sum_{\omega_0 , \cdots ,\omega_{-l}=1, \cdots , N} \left( \int_{B(\omega_0, \cdots , \omega_{-l})} \int_{\sigma_v(\omega')} \Psi d \nu_l(\omega') - \int \frac{1}{N^{l+1}}  \int_{\Delta^0_{\omega_0 \cdots \omega_{-l}}(t)} \Psi d \lambda(t) \right) \right| \\
&= \left| \sum_{\omega_0 , \cdots ,\omega_{-l}=1, \cdots , N} \int \int_{B(\omega_0, \cdots , \omega_{-l})} \left[ \int_{\sigma_v(\omega')} \Psi 
- \int_{\Delta^0_{\omega_0 \cdots \omega_{-l}}(t)} \Psi \right] d \nu_l(\omega')  d \lambda(t)  \right| \\
\end{split}
\end{equation*}

car $\lambda(\tau)=1$ et $\nu_l(B(\omega_0, \cdots , \omega_{-l}))= \frac{1}{N^{l+1}} $.

Pour $\omega'$ dans $B(\omega_0, \cdots , \omega_{-l})$, on a $\omega_0'=\omega_0 , \cdots , \omega_{-l}'= \omega_{-l}$, ce qui implique que $\sigma_v(\omega')$ et $\Delta^0_{\omega_0 \cdots \omega_{-l}}(t)$ sont dans
$\overline{\mathcal{A}^0_{ \omega_0 \cdots \omega_{-l}}}$. On a donc

$$d(\sigma_v(\omega'), \Delta^0_{\omega_0 \cdots \omega_{-l}}(t)) \leq 2h e^{- \gamma n(l+1) + 2 \gamma(l+1)}.$$

Si on applique le même raisonnement que pour la continuité de $\omega \longrightarrow \int_{\sigma_v(\omega)} \Psi$, on obtient que

$$ \left| \int_{\sigma_v(\omega')} \Psi 
- \int_{\Delta^0_{\omega_0 \cdots \omega_{-l}}(t)} \Psi \right| \leq \epsilon(l)$$

avec $\epsilon(l)$ indépendant de $t$ et $\omega'$ et qui converge vers $0$ quand $l$ tend vers l'infini.

Cela termine la démonstration du lemme.

\end{proof}

Reprenons la preuve du théorème \ref{general}. Il s'agit de montrer que $\int \chi {\tau_y}_* C_{\gamma}(\widehat{y})_* S_{\widehat{y}} \wedge \omega^{k-s}=0$.

En utilisant le lemme précédent, on a

\begin{equation*}
\begin{split}
&\int \chi {\tau_y}_* C_{\gamma}(\widehat{y})_* S_{\widehat{y}} \wedge \omega^{k-s}= \int S_{\widehat{y}} \wedge C_{\gamma}(\widehat{y})^* {\tau_y}^*(\chi \omega^{k-s})\\
&= \lim_{l \rightarrow \infty} \int \frac{1}{N^{l+1}} \sum_{\omega_0 , \cdots ,\omega_{-l}=1, \cdots , N} \int_{\Delta^0_{\omega_0 \cdots \omega_{-l}}(t)} C_{\gamma}(\widehat{y})^* {\tau_y}^*(\chi \omega^{k-s}) d \lambda(t)\\
&=\lim_{l \rightarrow \infty} \int \frac{1}{N^{l+1}} \sum_{\omega_0 , \cdots ,\omega_{-l}=1, \cdots , N} \int_{\tau_y(C_{\gamma}(\widehat{y})(\Delta^0_{\omega_0 \cdots \omega_{-l}}(t)))} \chi \omega^{k-s} d \lambda(t).\\
\end{split}
\end{equation*}

Les $\tau_y(C_{\gamma}(\widehat{y})(\Delta^0_{\omega_0 \cdots \omega_{-l}}(t)))$ sont disjoints et inclus par construction dans $f^{-nl}(\tau_y(C_{\gamma}(\widehat{y})(\Delta^0(t)))$. Cela implique que

\begin{equation*}
\begin{split}
\sum_{\omega_0 , \cdots ,\omega_{-l}=1, \cdots , N} \int [\tau_y(C_{\gamma}(\widehat{y})(\Delta^0_{\omega_0 \cdots \omega_{-l}}(t)))]  d \lambda(t) &\leq \int [f^{-nl}(\tau_y(C_{\gamma}(\widehat{y})(\Delta^0(t))))] d \lambda(t)\\
& \leq (f^{nl})^* {\tau_y}_* C_{\gamma}(\widehat{y})_* \beta^s.\\
\end{split}
\end{equation*}

Ainsi

\begin{equation*}
\begin{split}
&\int \frac{1}{N^{l+1}} \sum_{\omega_0 , \cdots ,\omega_{-l}=1, \cdots , N} \int_{\tau_y(C_{\gamma}(\widehat{y})(\Delta^0_{\omega_0 \cdots \omega_{-l}}(t)))} \chi \omega^{k-s} d \lambda(t)\\
&\leq \frac{1}{N^{l+1}} \int (f^{nl})^* {\tau_y}_* C_{\gamma}(\widehat{y})_* \beta^s \wedge \chi \omega^{k-s} \leq \frac{1}{N^{l+1}} \int C_1 (f^{nl})^* \omega^s \wedge \chi \omega^{k-s}\\
&=\frac{C_1 d_s^{nl}}{N^{l+1}}  \int \left( \frac{(f^{nl})^* \omega^s}{d_s^{nl}}  -T^+ \right) \wedge \chi \omega^{k-s} \leq  \frac{C_1 d_s^{nl}}{N^{l+1}} e^{-\alpha_2 nl}\\
\end{split}
\end{equation*}

pour $l$ assez grand, grâce à l'hypothèse de convergence exponentielle.

Cette dernière quantité est enfin égale à $ \frac{C_1 e^{nl \log d_s -\alpha_2 nl}}{e^{n(l+1) \log d_s - \rho n(l+1)}}$ qui converge bien vers $0$ quand $l$ tend vers l'infini.

On a donc montré que le support de ${\tau_y}_* C_{\gamma}(\widehat{y})_* S_{\widehat{y}}$ est inclus dans le support de $T^+$. Par le théorème \ref{codage}, cela implique que 

$$\mbox{dim}_{\Hcal}(\mbox{Support de } T^+ ) \geq 2(k-s) + \frac{h_{\mu}(f)- \rho}{\chi_1 + 4 \gamma}.$$

Si on fait tendre $\rho$ puis $\gamma$ vers $0$ (ce que l'on peut faire car le courant $T^+$ est indépendant de $\rho$ et $\gamma$), on obtient le théorème \ref{general}.

\section{\bf{Cas des applications birationnelles de Bedford-Diller}}

Dans ce paragraphe on reprend les résultats et les notations d'E. Bedford et J.  Diller (voir \cite{BD}).

Tout d'abord, ils supposent que $d_1 > 1$. Cela implique par un résultat de J. Diller et C. Favre (voir le théorème 5.1 de \cite{DF}) que $d_1$ est l'unique (en comptant avec multiplicité) valeur propre de module strictement plus grand que $1$ de l'action $f^*$ sur $H^{1,1}(X)$. On notera $\theta^+ \in H^{1,1}(X)$ un vecteur propre associé à la valeur propre $d_1$. On considère aussi $\theta^-$ l'analogue pour $f_*$ et on a $\langle \theta^+ , \theta^- \rangle > 0$ par  J. Diller et C. Favre (\cite{DF}). Dans la suite, on normalise $\theta^{\pm}$ et $\omega$ de sorte que

$$\langle \theta^+ , \theta^- \rangle = \langle \theta^+ , \omega \rangle = \langle \theta^- , \omega \rangle =1.$$

Soit $\omega_1, \cdots , \omega_N$ des formes kählériennes qui constituent une base de $H^{1,1}(X)$ et $\Omega$ l'espace constitué des combinaisons linéaires de ces $\omega_i$. On peut supposer pour simplifier que $\omega$ est l'un des $\omega_i$.

Si $\eta$ est un $(1,1)$ courant fermé, on peut écrire $\eta = \omega(\eta)+ dd^c p(\eta)$ où $\omega(\eta)$ est l'élément de $\Omega$ de même classe que $\eta$ dans $H^{1,1}(X)$.

Soit $\omega^+$ l'élément de $\Omega$ qui représente $\theta^+$ (c'est-à-dire $\omega^+=\omega(\theta^+)$).

Si $\eta$ est dans $\Omega$, on peut aussi écrire $\eta = c \omega^+ + \eta^{\perp}$ où $\eta^{\perp}$ est dans le sous-espace associé aux valeurs propres de $f^*$ autres que $d_1$. Comme on a $f_* \theta^- = d_1 \theta^-$, on en déduit que $c= \langle \eta, \theta^- \rangle$.

Par Bedford-Diller (voir les équations $(12)$ et $(13)$), si $\eta$ est une $(1,1)$ forme lisse fermée, on a

$$ {f^n}^* \eta = \omega({f^n}^* \eta) + d_1^n dd^c g_n^+(\eta)$$

avec

$$g_n^+(\eta)= \frac{1}{d_1^n} \left( \sum_{j=0}^{n-1} {f^{(n-1-j)}}^* \gamma^+({f^{j}}^* \eta ) + {f^{n}}^* p(\eta) \right).$$

Ici, on a posé $\gamma^+(\eta)= p(f^*(\omega(\eta)))$. Comme dans Bedford-Diller, on considère maintenant la fonction 

$$\gamma^+:= \gamma^+(\omega^+)=  p(f^*(\omega(\omega^+)))= p (f^*(\omega^+)).$$

Le courant de Green $T^+$ est égal à $T^+=\omega^+ + dd^c g^+$ avec $g^+= \sum_{j=0}^{\infty} \frac{\gamma^+ \circ f^j }{d_1^{j+1}}$.

Alors, on a

\begin{Pro}{\label{convergence}}

Soit $(X, \omega)$ une surface kählérienne compacte et $f: X \longrightarrow  X$ une application birationnelle qui vérifie la condition de Bedford-Diller. On note $d_1$ son premier degré dynamique et on suppose que $d_1 > 1$. Alors $\frac{(f^n)^* \omega}{d_1^n}$ converge exponentiellement vite vers le courant de Green $T^+$ dans le sens suivant: 

Soit $1 < \lambda < d_1$. Pour toute $(1,1)$-forme lisse $\Psi$, il existe $n_0$ tel que pour tout $n \geq n_0$ on ait:

$$\left| \langle  \frac{(f^n)^* \omega}{d_1^n} - T^+, \Psi \rangle \right| \leq \left( \frac{\lambda}{d_1} \right)^n.$$

\end{Pro}

\begin{proof}

La démonstration de cette proposition va découler de la proposition 2.4 et des idées de la preuve du théorème 2.5 de l'article d'E. Bedford et J. Diller (voir \cite{BD}).

Soit $1 < \lambda < d_1$ et $1 < \lambda' < \lambda$.

On décompose $\omega= c \omega^+ + \omega^{\perp}$. Comme $\langle \theta^- , \omega \rangle =1$, on a $c=1$.

Soit $\Psi$ une $(1,1)$-forme lisse.

Pour démontrer la proposition, on va montrer que $\langle \frac{{f^n}^* \omega^+}{d_1^n} - T^+ , \Psi \rangle$ et $\langle \frac{{f^n}^* \omega^{\perp}}{d_1^n}  , \Psi \rangle$ tendent vers $0$ exponentiellement vite. 

On a 

$${f^n}^* \omega^+ = \omega({f^n}^* \omega^+) + d_1^n dd^c g_n^+(\omega^+)= d_1^n \omega^+ + d_1^n dd^c g_n^+(\omega^+)$$

avec

$$g_n^+(\omega^+)= \frac{1}{d_1^n} \left( \sum_{j=0}^{n-1} {f^{(n-1-j)}}^* \gamma^+({f^{j}}^* \omega^+ ) + {f^{n}}^* p(\omega^+) \right).$$

Tout d'abord, $p(\omega^+)=0$ car $\omega^+ \in \Omega$. Ensuite,

$$\gamma^+({f^{j}}^* \omega^+ )=p(f^*(\omega({f^{j}}^* \omega^+)))=p(f^*(d_1^j \omega^+))=d_1^j \gamma^+.$$

On a ainsi,

$$g_n^+(\omega^+)= \frac{1}{d_1^n} \sum_{j=0}^{n-1} {f^{(n-1-j)}}^*(d_1^j \gamma^+  ) = \sum_{j=0}^{n-1} \frac{\gamma^+ \circ f^j}{d_1^{j+1}}.$$

D'où

$$I_1=\left| \langle T^+ - \frac{{f^n}^* \omega^+}{d_1^n}, \Psi \rangle \right|= \left| \langle dd^c \sum_{j \geq n} \frac{ \gamma^+ \circ f^j}{d_1^{j+1}} , \Psi \rangle \right|= \left| \langle \sum_{j \geq n} \frac{ \gamma^+ \circ f^j}{d_1^{j+1}} , dd^c \Psi \rangle \right|.$$

Maintenant, on peut trouver des constantes $C_1(\Psi)$ et $C_2(\Psi)$ telles que $0 \leq  dd^c \Psi + C_1(\Psi) \omega^2 \leq C_2(\Psi) \omega^2$.

En écrivant $dd^c \Psi = dd^c \Psi + C_1(\Psi) \omega^2 - C_1(\Psi) \omega^2$, on a donc

\begin{equation*}
\begin{split}
&\left| \langle \sum_{j \geq n} \frac{ \gamma^+ \circ f^j }{d_1^{j+1}} , dd^c \Psi \rangle \right| \leq \left| \langle \sum_{j \geq n} \frac{ \gamma^+ \circ f^j}{d_1^{j+1}} , dd^c \Psi + C_1(\Psi) \omega^2 \rangle \right| + \left|  \langle \sum_{j \geq n} \frac{ \gamma^+ \circ f^j }{d_1^{j+1}} , C_1(\Psi) \omega^2 \rangle       \right| \\
&\leq \sum_{j \geq n} \int \left| \frac{\gamma^+ \circ f^j }{d_1^{j+1}} \right| (dd^c \Psi + C_1(\Psi) \omega^2)  +  \sum_{j \geq n} \int \left| \frac{ \gamma^+ \circ f^j }{d_1^{j+1}} \right|C_1(\Psi) \omega^2 \\
&\leq \sum_{j \geq n} \int \left| \frac{\gamma^+ \circ f^j }{d_1^{j+1}} \right| C_2(\Psi) \omega^2  +  \sum_{j \geq n} \int \left| \frac{ \gamma^+ \circ f^j }{d_1^{j+1}} \right|C_1(\Psi) \omega^2 . \\
\end{split}
\end{equation*}

Par la proposition 2.4 de \cite{BD}, il existe une constante $C$ telle que pour tout $j \geq 0$ on ait $ \int | \gamma^+ \circ f^j | d \omega^2 \leq  C \lambda'^j.$

Ainsi,

$$I_1=\left| \langle T^+ - \frac{{f^n}^* \omega^+}{d_1^n}, \Psi \rangle \right| \leq \frac{ C_3(\Psi) \lambda'^n}{d_1^n}$$

pour tout $n \geq 0$.

Passons maintenant au terme $\langle \frac{{f^n}^* \omega^{\perp}}{d_1^n}  , \Psi \rangle$.

On a 

$$\frac{{f^n}^* \omega^{\perp}}{d_1^n} = \frac{\omega({f^n}^* \omega^{\perp})}{d_1^n} + dd^c g_n^+(\omega^{\perp})$$ 

avec

$$g_n^+(\omega^{\perp})= \frac{1}{d_1^n} \left( \sum_{j=0}^{n-1} {f^{(n-1-j)}}^* \gamma^+({f^{j}}^* \omega^{\perp} ) + {f^{n}}^* p(\omega^{\perp}) \right).$$

Comme $\omega^{\perp}$ est dans le sous-espace associé aux valeurs propres de $f^*$ autres que $d_1$ (et donc de modules inférieurs ou égaux à $1$), on a

$$\left| \langle \frac{\omega({f^n}^* \omega^{\perp})}{d_1^n}, \Psi \rangle \right| \leq C_4(\Psi) \frac{\lambda'^n}{d_1^n}$$

pour tout $n \geq 0$.

Pour finir, il reste à montrer que $| \langle dd^c g_n^+(\omega^{\perp}), \Psi \rangle |= | \langle g_n^+(\omega^{\perp}), dd^c \Psi \rangle |   $ converge exponentiellement vite vers $0$.

Remarquons tout d'abord que $p(\omega^{\perp})=0$ car $\omega^{\perp} = \omega - \omega^+ \in \Omega$. Ensuite, en écrivant comme précédemment $dd^c \Psi = dd^c \Psi + C_1(\Psi) \omega^2 - C_1(\Psi) \omega^2$, on a 

\begin{equation*}
\begin{split}
&| \langle g_n^+(\omega^{\perp}), dd^c \Psi \rangle |  \leq  \frac{1}{d_1^n} \sum_{j=0}^{n-1} \left| \int  {f^{(n-1-j)}}^* \gamma^+({f^{j}}^* \omega^{\perp} )  dd^c \Psi \right| \\
&\leq \frac{1}{d_1^n} \sum_{j=0}^{n-1} \left| \int  {f^{(n-1-j)}}^* \gamma^+({f^{j}}^* \omega^{\perp} )  (dd^c \Psi + C_1(\Psi) \omega^2) \right|\\
& + \frac{1}{d_1^n} \sum_{j=0}^{n-1} \left| \int  {f^{(n-1-j)}}^* \gamma^+({f^{j}}^* \omega^{\perp} )  C_1(\Psi) \omega^2 \right| \\
&\leq \frac{1}{d_1^n} \sum_{j=0}^{n-1} \int  {| f^{(n-1-j)}}^* \gamma^+({f^{j}}^* \omega^{\perp} )|  (dd^c \Psi + C_1(\Psi) \omega^2) \\
& + \frac{1}{d_1^n} \sum_{j=0}^{n-1}  \int  | {f^{(n-1-j)}}^* \gamma^+({f^{j}}^* \omega^{\perp} ) | C_1(\Psi) \omega^2  \\
&\leq \frac{1}{d_1^n} \sum_{j=0}^{n-1}  \int  |{f^{(n-1-j)}}^* \gamma^+({f^{j}}^* \omega^{\perp} )|  C_2(\Psi) \omega^2  + \frac{1}{d_1^n} \sum_{j=0}^{n-1}  \int | {f^{(n-1-j)}}^* \gamma^+({f^{j}}^* \omega^{\perp} ) |  C_1(\Psi) \omega^2.  \\
\end{split}
\end{equation*}

Comme on a 

$$\gamma^+(\eta)=p(f^*(\omega(\eta)))=p(f^*(\omega(\omega(\eta))))=\gamma^+(\omega(\eta)),$$

la proposition 2.4 de Bedford-Diller (voir \cite{BD}) donne que

$$| \langle g_n^+(\omega^{\perp}), dd^c \Psi \rangle | \leq \frac{1}{d_1^n} \sum_{j=0}^{n-1} C_5(\Psi) \lambda'^{n-1-j} \| \omega({f^{j}}^* \omega^{\perp}) \|.$$

Enfin, comme précédemment, il existe une constante $C(\omega^{\perp})$ telle que $\| \omega({f^{j}}^* \omega^{\perp}) \| \leq  C(\omega^{\perp}) \lambda'^j$ pour tout $j \geq 0$ et on obtient une majoration 

$$| \langle g_n^+(\omega^{\perp}), dd^c \Psi \rangle | \leq \frac{C_6(\Psi) n \lambda'^n}{d_1^n}.$$

Si on combine les résultats obtenus entre eux, on a donc

$$\left| \langle  \frac{(f^n)^* \omega}{d_1^n} - T^+, \Psi \rangle \right| \leq \frac{C_7(\Psi)n \lambda'^n}{d_1^n} $$

qui est bien plus petit que $ \left( \frac{\lambda}{d_1} \right)^n$ pour $n$ assez grand. Cela démontre la proposition.

\end{proof}

Grâce à cette proposition, on peut appliquer le théorème \ref{general} avec $k=2$, $s=1$, $\mu= T^+ \wedge T^-$ (qui est bien d'entropie maximale $\log d_1$ et qui intègre $\log d(x, I)$ par Bedford-Diller) et $\alpha_2= \log \left( \frac{d_1}{\lambda} \right)$. On obtient ainsi une minoration de $\mbox{dim}_{\Hcal} ( \mbox{Supp} T^+)$ par $2 + \frac{\log d_1}{\chi_1} > 2$. Cela démontre le théorème \ref{bedford-diller}.

\noindent Henry De Thélin, Université Paris 13, Sorbonne Paris Cité, LAGA, CNRS (UMR 7539), F-93430, Villetaneuse, France.  

\noindent Email: {\tt dethelin@math.univ-paris13.fr}

\end{document}